\documentclass[11pt,a4paper,leqno]{amsart}
\usepackage{xspace}
\usepackage{amsfonts,amssymb}
\usepackage{amsthm}
\usepackage[all]{xy}
\usepackage{textcomp}
\usepackage{mathrsfs}
\usepackage{stmaryrd}
\usepackage{bbm}
\usepackage{oldgerm}
\usepackage[T1]{fontenc}
\usepackage[latin1]{inputenc}
\theoremstyle{definition}

\newtheorem{rem}[subsection]{Remark}
\theoremstyle{plain}
\newtheorem{lemma}[subsection]{Lemma}

\newtheorem{prop}[subsection]{Proposition}
\newtheorem{thm}[subsection]{Theorem}

\theoremstyle{remark}

\numberwithin{equation}{subsection}

\newcommand{\beq}{\begin{equation}}
\newcommand{\eeq}{\end{equation}}

\newcommand{\ra}{\rightarrow}

\newcommand{\hra}{\hookrightarrow}

\newcommand{\xra}{\xrightarrow}

\newcommand{\Z}{\mathbf{Z}}

\newcommand{\Q}{\mathbf{Q}}

\newcommand{\D}{\mathbb{D}}

\newcommand{\bD}{\mathbf{D}}

\newcommand{\bA}{\mathbf{A}}
\newcommand{\fb}{\mathbf{f}}
\newcommand{\bg}{\mathbf{g}}
\newcommand{\bh}{\mathbf{h}}

\newcommand{\cA}{\mathscr{A}}
\newcommand{\cF}{\mathscr{F}}
\newcommand{\cG}{\mathcal{G}}

\newcommand{\cO}{\mathcal {O}}

\newcommand{\fX}{\mathscr{X}}
\newcommand{\fY}{\mathscr{Y}}
\newcommand{\Kb}{\overline{K}}
\newcommand{\Ks}{K^{\mathrm{sep}}}

\newcommand{\nb}{\underline{n}}

\newcommand {\resp}{\emph{resp.\xspace}}
\newcommand{\ie}{\emph{i.e.,} }

\newcommand{\FSA}{\mathbf{FEA}}
\newcommand{\Set}{\mathbf{Set}}
\newcommand{\m}{\mathfrak{m}}
\newcommand{\univ}{\mathrm{univ}}

\DeclareMathOperator{\Gal}{\mathrm{Gal}}

\DeclareMathOperator{\Lie}{\mathrm{Lie}}

\DeclareMathOperator{\End}{\mathrm{End}}

\DeclareMathOperator{\Spec}{\mathrm {Spec}}
\DeclareMathOperator{\Sp}{\mathrm{Sp}}

\begin{document}
\title[Bound of ramification filtration]{An upper bound on the Abbes-Saito filtration for finite flat group schemes and applications}
\author{Yichao TIAN}
\address{Mathematics Department, Fine Hall, Washington Road, Princeton, NJ, 08544, USA}
\email{yichaot@princeton.edu}
\date{May 16th, 2010}
\maketitle

\begin{abstract}
Let $\cO_K$ be a complete discrete valuation ring of residue characteristic $p>0$, and $G$ be a finite flat group scheme over $\cO_K$ of order a power of $p$. We prove in this paper that the Abbes-Saito filtration of $G$ is bounded by a linear function of the degree of $G$.  Assume $\cO_K$ has generic characteristic $0$ and the residue field of $\cO_K$ is perfect. Fargues constructed the higher level canonical subgroups for a ``not too supersingular'' Barsotti-Tate group $\cG$ over $\cO_K$. As an application of our bound, we prove that the canonical subgroup of $\cG$  of level $n\geq 2$ constructed by Fargues appears in the Abbes-Saito filtration of the $p^n$-torsion subgroup of $\cG$.

\end{abstract}
\vspace{1cm}

Let $\cO_K$ be a complete discrete valuation ring with  residue field $k$ of characteristic $p>0$ and   fraction field $K$. We denote by $v_{\pi}$ the valuation on $K$ normalized by $v_{\pi}(K^{\times})=\Z$.   Let $G$ be a finite and flat group scheme over $\cO_K$ of order a power of $p$ such that $G\otimes K$ is \'etale. We denote by $(G^{a}, a\in \Q_{\geq 0})$ the  Abbes-Saito filtration of $G$. This is  a decreasing and separated filtration of $G$ by finite and flat closed subgroup schemes. We refer the readers to \cite{AS02,AS,AM} for a full discussion, and to section 1 for a brief review of this filtration.    Let $\omega_{G}$ be the module of invariant differentials of $G$. The generic \'etaleness of $G$ implies that $\omega_G$ is a torsion $\cO_K$-module of finite type. There exist thus nonzero elements $a_1,\cdots, a_d\in \cO_{K}$ such that
\[\omega_{G}\simeq \bigoplus_{i=1}^d \cO_K/(a_i).\]
 We put $\deg(G)=\sum_{i=1}^dv_{\pi}(a_i)$, and call it the degree of $G$.  The aim of this note is to prove the following

\begin{thm}\label{main-thm}
Let $G$ be a finite and flat group scheme over $\cO_K$ of order a power of $p$ such that $G\otimes K$ is \'etale. Then we have $G^{a}=0$ for $a> \frac{p}{p-1}\deg(G)$.
\end{thm}

Our bound is quite optimal when $G$ is killed by $p$. Let $E_{\delta}=\Spec(\cO_K[X]/(X^p-\delta X))$ be the group scheme of Tate-Oort over $\cO_K$. We have $\deg(E_{\delta})=v_{\pi}(\delta)$, and an easy computation by Newton polygons gives \cite[Lemme 5]{Fa09}
$$E_\delta^{a}=\begin{cases}E_\delta &\text{if } 0\leq a\leq \frac{p}{p-1}\deg(E_{\delta})\\
0&\text{if } a>\frac{p}{p-1}\deg(E_{\delta}).
\end{cases}$$
However, our bound may be improved when $G$ is not killed by $p$ or $G$ contains many identical copies of a closed subgroup. In \cite[Thm. 7]{Hat}, Hattori proves that if $K$ has characteristic $0$ and $G$ is killed by $p^n$, then the Abbes-Saito filtration of $G$ is bounded by that of the multiplicative group $\mu_{p^n}$, \ie we have $G^a=0$ if $a>en+\frac{e}{p-1}$ where $e$ is the absolute ramification index of $K$. Compared with Hattori's result,  our bound has the advantage that it works in both characteristic $0$ and characteristic $p$, and that it is good if $\deg(G)$ is small.

The basic idea to prove \ref{main-thm} is to approximate general power series over $\cO_K$ by linear functions. First, we choose a  ``good'' presentation of the algebra of $G$ such that the defining equations of $G$ involve only terms of total degree $m(p-1)+1$ with $m\in \Z_{\geq 0}$ (Prop. \ref{prop-present}).  The existence of such a presentation is a consequence of the classical theory on $p$-typical curves of formal groups. With this good presentation, we can prove that the neutral connected component of the $a$-tubular neighborhood of $G$ is isomorphic to a closed rigid ball for $a> \frac{p}{p-1}\deg(G)$ (Lemma \ref{lemma-1}), and the only zero of the defining equations of $G$ in the neutral component is the unit section.

The motivation of our theorem comes from the theory of canonical subgroups. We assume  that $K$ has characteristic $0$, and the residue field $k$ is perfect of characteristic $p\geq 3$. Let $G$ be a Barsotti-Tate group of dimension $d\geq 1$ over $\cO_K$. If $G$ comes from an abelian scheme over $A$,  the canonical subgroup of level $1$ of $G$ was first constructed by Abbes and Mokrane in \cite{AM}. Then the author generalized their result to the  Barsotti-Tate case \cite{Ti}. We actually proved that if a Barsotti-Tate group $G$   over $\cO_K$  is ``not too supersingular'', a condition expressed explicitly as a bound on the Hodge height of $G$ (cf. \ref{Hdg-Height}), then a certain piece of the Abbes-Saito filtration of $G[p]$ lifts the kernel of Frobenius of the special fiber of $G$ \cite[Thm. 1.4]{Ti}. Later on, Fargues \cite{Fa09} gave another construction of the canonical subgroup of level $1$ by using Hodge-Tate maps, and his approach also allowed us to construct by induction the canonical subgroups of  level $n\geq 2$, \ie the canonical lifts of the kernel of $n$-th iteration of the Frobenius. He proved that the canonical subgroup of higher level appears in the Harder-Narasihman filtration of $G[p^n]$, which was introduced by him in \cite{Fa07}. It is conjectured that the canonical subgroup of higher level also appears in the Abbes-Saito filtration of $G[p^n]$. In this paper, we prove this conjecture  as a corollary of \ref{main-thm} (Thm. \ref{thm-can}). We use essentially the result of Fargues on the degree of the quotient of $G[p^n]$  by its canonical subgroup of level $n$ (see Thm. \ref{thm-Far}(i)).

\subsection{Notation} In this paper, $\cO_K$ will denote a complete discrete valuation ring with residue field $k$ of characteristic $p>0$, and with fraction field $K$. Let $\pi$ be a uniformizer of $\cO_K$, and $v_\pi$ be the valuation on $K$ normalized by $v_{\pi}(\pi)=1$.   Let $\Kb$ be an algebraic closure of $K$, $\Ks$ be the separable closure of $K$ contained in $\Kb$, and $\cG_K$ be the Galois group $\Gal(\Ks/K)$. We denote still by  $v_{\pi}$ the unique extension of the valuation to $\Kb$.

\section{Proof of Theorem \ref{main-thm}}

We recall first the definition of the filtration of Abbes-Saito for finite flat group schemes according to \cite{AM,AS}.

\subsection{} For a semi-local ring $R$, we denote by $\m_R$ its Jacobson radical. An algebra $R$ over $\cO_K$ is called \emph{formally of finite type}, if $R$ is semi-local, complete with respect to the $\m_R$-adic topology, Noetherian and $R/\m_R$ is finite over $k$. We say an $\cO_K$-algebra $R$ formally of finite type is formally smooth, if each of the factors of $R$ is formally smooth over $\cO_K$.

Let $\FSA_{\cO_K}$ be the category of finite, flat and generially \'etale  $\cO_K$-algebras, and $\Set_{\cG_K}$ be the category of finite sets endowed with a continuous action of the Galois group $\cG_K$. We have the  fiber functor
\[\cF:\FSA_{\cO_K}\ra \Set_{\cG_K},\]
which associates with an object $A$ of $\FSA_{\cO_K}$ the set $\Spec(A)(\Kb)$ equipped with the natural action of $\cG_K$. We define a filtration on the functor $\cF$ as follows. For each object $A$ in $\FSA_{\cO_K}$, we choose a presentation
\begin{equation}\label{present}0\ra I\ra \cA\ra A\ra 0,\end{equation}
where $\cA$ is an $\cO_K$-algebra formally of finite type and formally smooth. For any $a=\frac{m}{n}\in \Q_{>0}$ with $m$ prime to $n$, we define $\cA^a$ to be the $\pi$-adic completion of the subring $\cA[I^n/\pi^m]\subset \cA\otimes_{\cO_K} K$ generated  over $\cA$ by all the $f/\pi^m$ with $f\in I^n$. The  $\cO_K$-algebra $\cA^{a}$ is topologically of finite type, and the tensor product $\cA^a\otimes_{\cO_K} K$ is an affinoid algebra over $K$ \cite[Lemma 1.4]{AS}. We put $X^a=\Sp(\cA^{a}\otimes_{\cO_K}K)$, which is a smooth affinoid variety over $K$ \cite[Lemma 1.7]{AS}. We call it the \emph{$a$-th tubular neighborhood of  $\Spec(A)$  with respect to the presentation \eqref{present}}. The  $\cG_K$-set of  the geometric connected components of $X^a$, denoted by $\pi_0(X^a(A)_{\Kb})$, depends only on the $\cO_K$-algebra $A$ and the rational number $a$, but not on the  choice of the presentation \cite[Lemma 1.9.2]{AS}.  For rational numbers $b>a>0$, we have  natural inclusions of affinoid varieties $\Sp(A\otimes_{\cO_K}K)\hra X^b\hra X^a$, which induce natural morphisms $\Spec(A)(\Kb)\ra \pi_0(X^b(A)_{\Kb})\ra \pi_0(X^a(A)_{\Kb})$. For a morphism $A\ra B$ in $\FSA_{\cO_K}$, we can choose properly presentations of $A$ and $B$ so that we have a functorial map $\pi_0(X^a(B)_{\Kb})\ra \pi_0(X^a(A)_{\Kb})$.    Hence we get, for any $a\in \Q_{>0}$, a (contravariant) functor
\[\cF^a:\FSA_{\cO_K}\ra \Set_{\cG_K}\]
given by $A\mapsto \pi_0(X^a(A)_{\Kb})$. We have natural morphisms of functors $\phi_a:\cF \ra \cF^a$, and $\phi_{a,b}:\cF^b\ra \cF^a$ for rational numbers $b>a>0$ with $\phi_a=\phi_{b,a}\circ\phi_b$. For any $A$ in $\FSA_{\cO_K}$, we have $\cF(A)\xra{\sim}\varprojlim_{a\in \Q_{>0}}\cF^a(A)$ \cite[6.4]{AS02};  if $A$ is a complete intersection over $\cO_K$, the map $\cF(A)\ra \cF^a(A)$ is surjective for any $a$ \cite[6.2]{AS02}.

\subsection{} Let $G=\Spec(A)$ be a finite and flat group scheme over $\cO_K$ such that $G\otimes K$ is \'etale over $K$, and $a\in \Q_{>0}$. The group structure of $G$ induces a group structure on $\cF^a(A)$, and the natural map $G(\Kb)=\cF(A)\ra \cF^{a}(A)$ is a homomorphism of groups. Hence the kernel $G^a(\Kb)$ of $G(\Kb)\ra \cF^a(A)$
 is a $\cG_K$-invariant subgroup of $G(\Kb)$, and it defines a closed subgroup scheme $G^a_K$ of the generic fiber $G\otimes K$. The scheme theoretic closure of $G^a_K$ in $G$, denoted by $G^a$, is a closed subgroup of $G$ finite and flat over $\cO_K$. Putting $G^0=G$,  we get a decreasing and separated filtration $(G^a,a\in \Q_{\geq 0})$ of $G$ by finite and flat closed subgroup schemes. We call it \emph{Abbes-Saito filtration} of $G$. For any real number $a\geq 0$, we put $G^{a+}=\cup_{b\in\Q_{>a}}G^a$.

Assume $G$ is connected, \ie the ring $A$ is local. Let
\begin{equation}\label{formal-present}0\ra I\ra \cO_K[[X_1,\cdots, X_d]]\ra A\ra 0\end{equation}
 be a presentation of $A$ by  the ring of formal power series
such that the unit section of $G$ corresponds to the point $(X_1,\cdots, X_d)=(0,\cdots, 0)$. Since $A$ is a relative complete intersection over $\cO_K$,  $I$ is generated by $d$ elements $f_1,\cdots, f_d$.   For $a\in \Q_{>0}$, the $\Kb$-valued points of the $a$-th tubular neighborhood  of $G$ are given by
\begin{equation}\label{tubular}
X^a(\Kb)=\bigl\{(x_1,\cdots,x_d)\in \m_{\Kb}^d\;|\; v_{\pi}(f_i(x_1,\cdots,x_d))\geq a \;\text{for }1\leq i\leq d\bigr\},\end{equation}
where $\m_{\Kb}$ is the maximal ideal of $\cO_{\Kb}$.
The subset $G(\Kb)\subset X^a(\Kb)$ corresponds to the zeros of the $f_i$'s.
Let $X^a_0$ be the connected component of $X^a$ containing $0$. Then the subgroup $G^a(\Kb)$ is the intersection of $X^a_0(\Kb)$ with $G(\Kb)$.

The basic properties of Abbes-Saito filtration that we need are summarized as follows.

\begin{prop}[\cite{AM} 2.3.2, 2.3.5]\label{prop-AS} Let $G$ and $H$ be finite and flat group schemes, generically \'etale over $\cO_K$, $f:G\ra H$ be a homomorphism of group schemes.

\emph{(i)} $G^{0+}$ is the connected component of $G$, and we have $(G^{0+})^{a}=G^a$ for any $a\in \Q_{> 0}$.

\emph{(ii)} For $a\in \Q_{>0}$, $f$ induces a canonical homomorphism $f^a:G^a\ra H^a$. If $f$ is flat and surjective, then $f^a(\Kb):G^a(\Kb)\ra H^{a}(\Kb)$ is surjective.
\end{prop}

Now we return to the proof of Theorem \ref{main-thm}.

\begin{lemma}
Let $R$ be a $\Z_p$-algebra, $\fX$ be a formal group of dimension $d$ over $R$. Then

\emph{(i)}  the ring $\Z_p$ acts naturally on $\fX$, and its image in $\End_R(\fX)$ lies in the center of $\End_{R}(\fX)$;

\emph{(ii)}  there exist parameters  $(X_1,\cdots, X_d)$ of $\fX$, such that we have $[\zeta](X_1,\cdots, X_d)=(\zeta X_1,\cdots \zeta X_d)$ for any $(p-1)$-th root of unity $\zeta\in \Z_p$.
\end{lemma}

\begin{proof} This is actually a classical result on formal groups. In the terminology of \cite{Ha}, $\fX$ is necessarily isomorphic to a $p$-typical formal group over $R$ \cite[16.4.14]{Ha}. This means that  $\fX$ is deduced by base change from the universal $p$-typical formal group $\fX^{\mathrm{univ}}$ (denoted by $F_V(X,Y)$ in \cite[15.2.8]{Ha}) over $\Z_p[V]=\Z_p[V_{i}(j,k); i\in \Z_{\geq 0}, j,k=1,\cdots,d]$, where the $V_i(j,k)$'s are free variables. So we are reduced to proving the Lemma for $\fX^{\univ}$. If $X$ and $Y$ are short for the column vectors $(X_1,\cdots, X_d)$ and $(Y_1,\cdots,Y_d)$ respectively,  the formal group law on $\fX^{\univ}$ is determined by
$$F_V(X,Y)=f_{V}^{-1}(f_V(X)+f_V(Y)),\quad \text{with}\; f_V(X)=\sum_{i=0}^{\infty}a_{i}(V)X^{p^i},$$
where $a_i(V)$'s are certain $d\times d$ matrices with coefficients in $\Q_p[V]$ with $a_1(V)$ invertible, $X^{p^i}$ is short for $(X_1^{p^i},\cdots, X_d^{p^i})$, and $f^{-1}_V$ is the unique $d$-tuple of power series in $(X_1,\cdots, X_d)$ with coefficients in $\Q_p[V]$ such that $f^{-1}_V\circ f_V=1$ \cite[10.4]{Ha}. We note that $F_V(X,Y)$ is a $d$-tuple of power series with coefficient in $\Z_p[V]$, although $f_V(X)$ has coefficients in $\Q_p[V]$ \cite[10.2(i)]{Ha}. Via approximation by integers, we see easily that the multiplication by an element $\xi\in \Z_p$ can be well defined as $[\xi](X)=f_{V}^{-1}(\xi f_V(X))$. This proves (i). Statement (ii) is an immediate consequence of the fact that $f_V(X)$ involves just $p$-powers of $X$.
\end{proof}

\begin{prop}\label{prop-present}
Let $G=\Spec(A)$ be a connected finite and flat group scheme over $\cO_K$ of order a power of $p$. Then there exists a presentation of $A$ of type \eqref{formal-present} such that the defining equations $f_i$ for $1\leq i\leq d$ have the form
\[f_i(X_1,\cdots,X_d)=\sum_{|n|\geq 1}^{\infty} a_{i,\nb}X^{\nb}\quad \quad \text{with $a_{i,\nb}=0\;$ if $(p-1)\nmid (|\nb|-1)$,}\]
where $\nb=(n_1,\cdots,n_d)\in (\Z_{\geq 0})^d$ are multi-indexes, $|\nb|=\sum_{j=1}^dn_j$, and $X^{\nb}$ is short for $\prod_{j=1}^dX_j^{n_j}.$
\end{prop}

\begin{proof} By a theorem of Raynaud \cite[3.1.1]{BBM}, there is a projective abelian variety $V$ over $\cO_K$, and an embedding of group schemes $j: G\hra V$. Let $\fX$ be the formal completion of $V$ along its unit section. This is a formal group over $\cO_K$. Since $G$ is connected,  then $j$ induces an embedding $i:G\hra \fX$. We denote by $\fY$ the quotient of $\fX$ by $G$, and by $\phi: \fX\ra \fY$  the canonical isogeny. Let $(X_1,\cdots, X_d)$ (\resp ($Y_1,\cdots, Y_d$)) be parameters  of $\fX$ (\resp $\fY$) satisfying the lemma above. The isogeny $\phi$ is thus given by
\[(X_1,\cdots, X_d)\mapsto (f_1(X_1,\cdots,X_d),\cdots, f_d(X_1,\cdots, X_d)),\]  where $f_i=\sum_{|\nb|\geq 1}a_{i,\nb}X^{\nb}\in \cO_K[[X_1,\cdots, X_d]]$. Since for any $(p-1)$-th root of unity $\zeta\in \Z_p$ we have $f_i(\zeta X_1,\cdots, \zeta X_d)=\zeta f_i(X_1,\cdots, X_d)$, it's easy to see that $a_{i,\nb}=0$ if $(p-1)\nmid (|\nb|-1)$.

\end{proof}

\subsection{Proof of Theorem \ref{main-thm}} Let $H=G^{0+}$ be the connected component of $G$. By \ref{prop-AS}(i), we have $G^a=H^a$ for $a\in \Q_{>0}$. On the other hand,  from the exact sequence of group schemes $0\ra H\ra G\ra G/H\ra 0$, it follows that the sequence of finite $\cO_K$-modules
$$0\ra \omega_{G/H}\ra \omega_{G}\ra \omega_H\ra 0$$ is exact. Since $G/H$ is \'etale,  we have $\omega_{G/H}=0$ and hence $\deg(G)=\deg(H)$. Up to replacing $G$ by $H$, we may assume that $G=\Spec(A)$ is connected.

We choose a presentation of $A$ as in Prop. \ref{prop-present} so that we have an isomorphism of $\cO_K$-algebras
\[A\simeq \cO_K[[X_1,\cdots,X_d]]/(f_1,\cdots,f_d)\]
where \[f_i(X_1,\cdots, X_d)=\sum_{j=1}^da_{i,j}X_j+\sum_{|\nb|\geq p}a_{i,\nb}X^{\nb}.\]
 Then we have
$$\Omega^1_{A/\cO_K}\simeq \bigl(\bigoplus_{i=1}^d A dX_i\bigr)/(df_1,\cdots,df_d).$$
Since $\omega_G\simeq e^*(\Omega^1_{A/\cO_K})$, where $e$ is the unit section of $G$, we get
$$\omega_G\simeq \bigl(\oplus_{i=1}^d\cO_K dX_i\bigr)/(\sum_{1\leq j\leq d}a_{i,j}dX_j)_{1\leq i\leq d}.$$
 In particular, if $U$  denotes the matrix $(a_{i,j})_{1\leq i,j\leq d}$, then we have $\deg(G)=v_{\pi}(\det(U))$.

 For any rational number $\lambda$, we denote by $\bD^d(0,|\pi|^\lambda)$ (\resp $\D^d(0,|\pi|^\lambda)$) the rigid analytic closed (\resp open) disk of dimension $d$ over $K$ consisting of points $(x_1,\cdots,x_d)$ with $v_{\pi}(x_i)\geq \lambda$ (\resp $v_{\pi}(x_i)> \lambda$) for $1\leq i\leq d$; we put $\bD^d(0,1)=\bD^d(0,|\pi|^0)$ and $\D^d(0,1)=\D^d(0,|\pi|^0)$. Let $a>\frac{p}{p-1}\deg(G)$ be a rational number,   $X^a$ be the $a$-th tubular neighborhood of $G$ with respect to the chosen presentation.  By \eqref{tubular}, we have a cartesian diagram of rigid analytic spaces
 \begin{equation}\label{diag-tubular}
 \xymatrix{X^a\ar@{^(->}[r]\ar[d]^{\fb}&\D^d(0,1)\ar[d]^{\fb=(f_1,\cdots,f_d)}\\
 \bD^d(0,|\pi|^{a})\ar@{^(->}[r]&\D^d(0,1),}
 \end{equation} where horizontal arrows are  inclusions, and $\fb(y_1,\cdots,y_d)=(f_1(y_1,\cdots,y_d),\cdots,f_d(y_1,\cdots, y_d))$. Let $X^a_0$ be the connected component of $X^a$ containing $0$. By the discussion below \eqref{tubular},  we just need to  prove that  $0$ is the only zero of the $f_i$'s contained in $X^a_0$.

 Let $V=(b_{i,j})_{1\leq i,j\leq d}$ be  the unique $d\times d$ matrix with coefficients in $\cO_K$ such that $UV=VU=\det(U)I_d$, where $I_d$ is the $d\times d$ identity matrix. If $\bA^d_K$ denotes the $d$-dimensional rigid affine space over $K$, then $V$ defines an isomorphism of rigid spaces
 $$\bg: \bA^d_K\ra \bA^d_K;  \quad\quad (x_1,\cdots, x_d)\mapsto (\sum_{j=1}^db_{1,j}x_j,\cdots, \sum_{j=1}^db_{d,j}x_j).$$
 It's clear that $\bg(\D^{d}(0,1))\subset \D^{d}(0,1)$, so that $\fb$ is defined on $\bg(\D^d(0,1))$.
The composite morphism $\fb\circ \bg: \D^d(0,1)\ra \D^d(0,1)$ is given by
\begin{equation}\label{formula-1}
(x_1,\cdots,x_d)\mapsto (\det(U)x_1+R_1,\cdots, \det(U)x_d+R_d),
\end{equation}
where  $R_i=\sum_{|\nb|\geq p}a_{i, \nb}\prod_{j=1}^d(\sum_{k=1}^db_{j,k}x_k)^{n_j}$  involves only terms of order $\geq p$ for  $1\leq i\leq d$.  For $1\leq i\leq d$, we have   basic estimations
\begin{equation}\label{basic-estim}
v_{\pi}(\det(U)x_i)=\deg(G)+v_{\pi}(x_i)\quad \text{and}\quad v_{\pi}(R_i)\geq p\min_{1\leq j\leq d}\{v_{\pi}(x_j)\}.
\end{equation}

\begin{lemma}\label{lemma-1}
For any rational number $a> \frac{p}{p-1}\deg(G)$,
the map $\bg$ induces an isomorphism of affinoid rigid spaces
\[\bg: \bD^d(0, |\pi|^{a-\deg(G)})\xra{\sim} X^a_0.\]
\end{lemma}

Assuming this Lemma for a moment, we can complete the proof of \ref{main-thm} as follows. Consider the composite
\[\bh=\fb\circ\bg|_{\bD^d(0,|\pi|^{a-\deg(G)})}: \bD^d(0,|\pi|^{a-\deg(G)})\xra{\sim}X^a_0\hra X^a\xra{\fb} \bD^d(0,|\pi|^{a}).\]
In order to complete the proof of \ref{main-thm}, we just need to prove that the inverse image $\bh^{-1}(0)=\{0\}$. Let $(x_1,\cdots,x_d)$ be a point of $\bD^d(0,|\pi|^{a-\deg(G)})$, and $(z_1,\cdots,z_d)=\bh(x_1,\cdots,x_d)$. We may assume $v_{\pi}(x_1)=\min_{1\leq i\leq d}\{v_{\pi}(x_i)\} $. We have $v_{\pi}(x_1)\geq a-\deg(G)>\frac{1}{p-1}\deg(G)$ by the  assumption on $a$. It follows thus from
\eqref{basic-estim} that
$$\quad v_{\pi}(R_1)\geq pv_{\pi}(x_1)> \deg(G)+v_{\pi}(x_1)=v_{\pi}(\det(U)x_1).$$
Hence, we deduce from \eqref{formula-1} that
$v_{\pi}(z_1)=\deg(G)+v_{\pi}(x_1)$. In particular, $z_1=0$ if and only if $x_1=0$. Therefore, we have $\bh^{-1}(0)=\{0\}$. This achieves the proof of Theorem \ref{main-thm}.

\begin{proof}[Proof of \ref{lemma-1}]
Let $\epsilon$ be any rational number with $0<\epsilon <\frac{p-1}{p}a-\deg(G)$. We will prove that
\[\bD^d(0,|\pi|^{a-\deg(G)})=\bD^d(0,|\pi|^{a-\deg(G)-\epsilon})\cap \bg^{-1}(X^a).\]
This will imply that $\bD^d(0,|\pi|^{a-\deg(G)})$ is a connected component of $\bg^{-1}(X^a)$. Since $\bg:\bA^d_K\ra \bA^d_K$ is an isomorphism, the lemma will follow immediately.

We prove first the inclusion ``$\subset$''. It suffices to show $\bg(\bD^d(0,|\pi|^{a-\deg(G)}))\subset X^a$. Let $(x_1,\cdots, x_d)$ be a   point of $\bD^d(0, |\pi|^{a-\deg(G)})$. By \eqref{diag-tubular}, we have to check that $(z_1, \cdots,z_d )=\fb(\bg(x_1,\cdots, x_d))$ lies in $\bD^d(0,|\pi|^a)$. We get from \eqref{basic-estim} that $v_{\pi}(\det(U)x_i)=\deg(G)+v_{\pi}(x_i)\geq a$ and $v_{\pi}(R_i)\geq p(a-\deg(G))$. As
 $a> \frac{p}{p-1}\deg(G)$, we have  $v_{\pi}(R_i)> a$. It follows from \eqref{formula-1} that
$$v_{\pi}(z_i)\geq \min\{v_{\pi}(\det(U)x_i,v_{\pi}(R_i)\}\geq  a.$$
This proves $(z_1,\cdots, z_d)$ is contained in $\bD^d(0,|\pi|^a)$, hence we have  $\bg(\bD^d(0,|\pi|^{a-\deg(G)}))\subset X^a$.

 To prove the inclusion ``$\supset$'', we just need to verify that every point in $\bD^d(0,|\pi|^{a-\deg(G)-\epsilon})$ but outside $\bD^d(0,|\pi|^{a-\deg(G)})$ does not lie in $\bg^{-1}(X^a)$. Let $(x_1,\cdots, x_d)$ be such a point. We may assume that
 \begin{equation}\label{estimation-1}
 a-\deg(G)-\epsilon\leq v_{\pi}(x_1)<a-\deg(G)\quad \text{and}\quad v_{\pi}(x_i)\geq a-\deg(G)-\epsilon\;\text{ for }2\leq i\leq d.
 \end{equation}
 Let $(z_1,\cdots, z_d)=(\det(U)x_1+R_d,\cdots, \det(U)x_d+R_d)$ be the image of $(x_1,\cdots,x_d)$ under the composite $\fb\circ \bg$. According to \eqref{diag-tubular}, the proof will be completed if we can prove that $(z_1,\cdots,z_d)$ is not in $\bD^d(0,|\pi|^{a})$.  From \eqref{basic-estim} and \eqref{estimation-1}, we get
$v_{\pi}(\det(U)x_1)=\deg(G)+v_{\pi}(x_1)<a$
and $v_{\pi}(R_1)\geq p(a-\deg(G)-\epsilon)$. Thanks to the assumption on $\epsilon$, we have $p(a-\deg(G)-\epsilon)>a$, so $v_{\pi}(z_1)=v_{\pi}(\det(U)x_1)<a$. This shows that $(z_1,\cdots, z_d)$ is not in $\bg^{-1}(X^a)$,  hence the proof of the lemma is complete.

\end{proof}

\section{Applications to Canonical subgroups}
In this section, we suppose the fraction field $K$ has characteristic $0$ and the residue field $k$ is perfect of characteristic $p\geq 3$. Let $e$ be the absolute ramification index of $\cO_K$. For any rational number $\epsilon >0$, we denote by $\cO_{K,\epsilon}$ the quotient of $\cO_K$ by the ideal consisting of elements with $p$-adic valuation greater or equal than $ \epsilon$.

\subsection{}\label{Hdg-Height} First we recall some results on the canonical subgroups according to \cite{AM}, \cite{Ti} and \cite{Fa09}. Let $v_p:\cO_K/p\ra [0,1]$ be the truncated $p$-adic valuation (with $v_p(0)=1$). Let $G$ be a truncated Barsotti-Tate group of level $n\geq 1$ non-\'etale over $\cO_K$, $G_1=G\otimes_{\cO_K}(\cO_{K}/p)$.  The Lie algebra of $G_1$, denoted by $\Lie(G_1)$ is a finite free $\cO_K/p$-module. The Verschiebung homomorphism $V_{G_1}:G^{(p)}_1\ra G_1$ induces a semi-linear endomorphism $\varphi_{G_1}$ of $\Lie(G_1)$. We choose a basis of $\Lie(G_1)$ over $\cO_K/p$, and let $U$ be the matrix of $\varphi$ under this basis. We define the Hodge height of $G$, denoted by $h(G)$, to be the truncated $p$-adic valuation of $\det(U)$. We note that the definition of $h(G)$ does not depend on the choice of $U$. The Hodge height of $G$ is an analog of the Hasse invariant in mixed characteristic, and we have $h(G)=0$ if and only if $G$ is ordinary.

\begin{thm}[\cite{Fa09} Th\'eo. 4]\label{thm-Far-1} Let $G$ be a truncated Barsotti-Tate group of level $1$ over $\cO_K$ of dimension $d\geq 1$ and height $h$. Assume $h(G)< \frac{1}{2}$ if $p\geq 5$ and $h(G)<1/3$ if $p=3$.

\emph{(i)} For any rational number $\frac{ep }{p-1}h(G)<a\leq \frac{ep}{p-1}(1-h(G))$, the finite flat subgroup $G^{a}$ of $G$ given by the Abbes-Saito filtration has rank $p^d$.

\emph{(ii)} Let $C$ be the subgroup $G^{\frac{ep}{p-1}(1-h(G))}$ of $G$.  We have $\deg(G/C)=e h(G)$.

\emph{(iii)} The subgroup $C\otimes\cO_{K,1-h(G)}$ coincides with the kernel of the Frobenius homomorphism of $G\otimes\cO_{K,1-h(G)}$. Moreover, for any rational number $\epsilon$ with $\frac{h(G)}{p-1}<\epsilon\leq 1-h(G)$, if $H$ is a finite and flat closed subgroup of $G$ such that $H\otimes \cO_{K,\epsilon}$ coincides with the kernel of Frobenius of $G\otimes \cO_{K,\epsilon}$, then we have  $H=C$.
\end{thm}
The subgroup $C$ in this theorem, when it exists, is called the \emph{canonical subgroup (of level 1)} of $G$.

\begin{rem}
(i) The conventions here  are slightly different from those in \cite{Fa09}. The Hodge height is called Hasse invariant in \emph{loc. cit.}, while we choose to follow the terminologies in \cite{AM} and \cite{Ti}. Our index of Abbes-Saito filtration and the degree of $G$  are $e$ times those in \cite{Fa09}.

(ii) Statement (iii) of the theorem is not explicitly stated in \cite[Th\'eo. 4]{Fa09}, but it's an easy consequence of \emph{loc. cit.} Prop. 11.
\end{rem}

For the canonical subgroups of higher level, we have
\begin{thm}[\cite{Fa09} Th\'eo. 6]\label{thm-Far} Let $G$ be a truncated Barsotti-Tate group of level $n$ over $\cO_K$ of dimension $d\geq 1$ and height $h$. Assume $h(G)<\frac{1}{3^n}$ if $p=3$ and $h(G)<\frac{1}{2p^{n-1}}$ if $p\geq 5$.

\emph{(i)} There exists a unique closed subgroup of $G$ that is finite and flat over $\cO_K$ and satisfies
\begin{itemize}
\item $C_n(\Kb)$ is free of rank $d$ over $\Z/p^n\Z$.

\item For each integer $i$ with $1\leq i\leq n$, let $C_i$ be the scheme theoretic closure of $C_n(\Kb)[p^i]$ in $G$.  Then the subgroup $C_i\otimes \cO_{K,1-p^{i-1}h(G)}$ coincides with the kernel of the $i$-th iterated Frobenius of $G\otimes \cO_{K,1-p^{i-1}h(G)}$.
\end{itemize}

\emph{(ii)} We have $\deg(G/C_n)=\frac{e(p^n-1)}{p-1}h(G)$.
\end{thm}

The subgroup $C_n$ in the theorem above is called the canonical subgroup of level $n$ of $G$.  Fargues actually  proves that  $C_n$ is  a certain piece of the Harder-Narasimhan filtration of $G$.  The aim of this section is to show that $C_n$ appears also in the Abbes-Saito filtration.

\begin{thm}\label{thm-can}
Let $G$ be a truncated Barsotti-Tate group of level $n$ over $\cO_K$ satisfying the assumptions in $\ref{thm-Far}$, and $C_n$ be its canonical subgroup of level $n$. Then for any rational number $a$ satisfying $\frac{ep(p^n-1)}{(p-1)^2}h(G)<a\leq \frac{ep}{p-1}(1-h(G))$, we have $G^{a}=C_n$.
\end{thm}
\begin{proof}
We proceed by induction on $n$. If $n=1$, the theorem is \ref{thm-Far-1}(i). We suppose $n\geq 2$ and the theorem is valid for truncated Barsotti-Tate groups of level $n-1$. For each integer  $i$ with $1\leq i\leq n$, let $G_i$ denote the scheme theoretic closure of $G(\Kb)[p^i]$ in $G$,  and $C_i$ the scheme theoretic closure of $C_n(\Kb)[p^i]$ in $C_n$. By Theorem \ref{thm-Far}(i), it's clear that $C_i$ is the canonical subgroup of level $i$ of $G_i$.  Let $a$ be a rational number with $\frac{ep(p^n-1)}{(p-1)^2}h(G)<a\leq \frac{ep}{p-1}(1-h(G))$. By the induction hypothesis and the functoriality of Abbes-Saito filtration \ref{prop-AS}(ii), we have $C_{n-1}(\Kb)=G_{n-1}^a(\Kb)\subset G^a(\Kb)$, and  the image of $G^{a}(\Kb)$ in $G_1(\Kb)$ is exactly $C_1(\Kb)=G_1^a(\Kb)$. Note that we have a commutative diagram of exact sequences of groups
\[\xymatrix{0\ar[r]&C_{n-1}(\Kb)\ar[r]\ar@{^(->}[d]&C_n(\Kb)\ar[r]\ar@{^(->}[d]&C_1(\Kb)\ar[r]\ar@{^(->}[d]&0\\
0\ar[r]&G_{n-1}(\Kb)\ar[r]&G(\Kb)\ar[r]^{\times p^{n-1}}&G_1(\Kb)\ar[r]&0,}\]
where vertical arrows are natural inclusions. So we have $C_n(\Kb)\subset G^a(\Kb)$. On the other hand,  Theorems \ref{main-thm} and \ref{thm-Far}(ii) imply that $(G/C_n)^a(\Kb)=0$ as $a>\frac{ep(p^n-1)}{(p-1)^2}h(G)=\frac{p}{p-1}\deg(G/C_n)$. Therefore, we get $G^a(\Kb)\subset C_n(\Kb)$ by \ref{prop-AS}(ii). This completes the proof.

\end{proof}

\end{document}